  \documentclass[]{IEEEtran}
 
\usepackage[T1]{fontenc}
\usepackage[utf8]{inputenc}
 
\usepackage{cite}    
\usepackage{hyperref}      

\usepackage[dvips]{graphicx}  

\usepackage{amsmath,amsfonts,amssymb}

\usepackage{color}
\definecolor{mydarkblue}{rgb}{0,0.0,0.85}
\newcommand\bl[1]{#1}

\DeclareMathOperator*{\minimize}{minimize}

\begin{document}

\title{Tikhonov Regularization of Circle-Valued Signals}

\author{\begin{minipage}{\textwidth}\begin{center}{\Large{Laurent~Condat}}\\[10pt]  \Large
Author's final version.
Published in IEEE Transactions on Signal Processing, June 2022. \url{https://doi.org/10.1109/TSP.2022.3179816}
 \end{center}\end{minipage}
\thanks{L. Condat is with the Visual Computing Center, King Abdullah University of Science and Technology (KAUST), Thuwal 23955-6900, Kingdom of Saudi Arabia.
Contact: see https://lcondat.github.io/}
}

\maketitle
\begin{abstract}
It is common to have to process signals or images whose values are cyclic and can be represented as points on the complex circle, like wrapped phases, angles, orientations, or color hues.  
We consider a Tikhonov-type regularization model to smoothen or interpolate circle-valued signals defined on arbitrary graphs. We propose a convex relaxation of this nonconvex problem as a semidefinite program, and an efficient algorithm to solve it.
\end{abstract}

\section{Introduction}

\IEEEPARstart{I}{n} a wide range of applications, one has to deal with signals or images with cyclic, or circular, values, like phases, angles, orientations, or color hues, which are defined modulo $\pi$ or $2\pi$.  Cyclic data appear, for instance, in interferometric synthetic aperture radar \cite{ros00}, color image restoration in HSV or LCh spaces, profilometry~\cite{mos15}, Magnetic Resonance Imaging \cite{lan08}, biology, with data on the bacterial flagellar motor \cite{sow05}, in times series of wind directions \cite{dav02}, or in social sciences \cite{cre18}.

A cyclic value can be represented by a point on the complex circle; that is, a complex number 
%with amplitude 1, 
of the form $e^{j\omega}$, for some phase $\omega\in\mathbb{R}$, where $j=\sqrt{-1}$. Equivalently, the value is represented by its wrapped phase $\omega \in(-\pi,\pi]$, and the signal presents artificial $2\pi$ jumps when the values cross the $\pi$ or $-\pi$ boundaries. % of the interval $(-\pi,\pi]$.
Thus, to denoise or estimate circle-valued data, an option is to \emph{unwrap} the phase map to remove these artificial discontinuities, by estimating the  lost integer multiples of $2\pi$ in the phase values. Then the unwrapped signal or image can be processed using standard techniques for scalar data. Unfortunately, image unwrapping is a notoriously difficult %, often NP-hard, 
problem~\cite{
%mac83,%hun93,
%cus95,sal97,ghi98,
yin06,bio072,con19}, and the unwrapping process is prone to errors, so that it is preferable to process circle-valued data by keeping them on the circle.

We consider the general setting, where a signal is defined on a graph, with values located at the nodes. Two values are adjacent if there is an edge between their nodes. A 2-D image is a particular case with edges between every pair of neighboring pixels horizontally and vertically, forming a square grid. Then, to regularize signals on graphs, it is natural to promote the property that adjacent values are close to each other, in some sense. For scalar values, Tikhonov regularization consists in penalizing the squared differences of adjacent values \bl{and total variation (TV) regularization~\cite{cha101,con132,con17} consists, instead, in penalizing the absolute values of these  differences. In this work, we focus on Tikhonov regularization for circle-valued signals.  
A few methods have been proposed for TV  regularization of circle-valued signals  \cite{cre13,wei14,sto161}. There seems to be no available method for Tikhonov regularization of circle-valued data, with the exception of \cite{wei14}, where an iterative method based on the proximal point algorithm is proposed; it converges to the global solution for this type of problems on Hadamard manifolds, which the circle is not. Thus, we tackle the difficult nonconvex problem of Tikhonov regularization for circle-valued signals on graphs, by proposing a new convex relaxation.}

The paper is organized as follows: in Section~\ref{sec2}, we discuss different formulations for the regularization of circle-valued signals. In Section~\ref{secpro}, we propose a convex relaxation of the considered nonconvex problem, and in Section~\ref{sec4}, we propose an algorithm to solve it. In Section~\ref{sec5}, we illustrate the benefits of the proposed approach with several experiments.

\section{Tikhonov Smoothing for Circle-Valued Signals}\label{sec2}

\subsection{Circle-Valued Signals on Graphs}

Let $\mathbb{S}=\{z\in\mathbb{C}\ :\ |z|=1\}$ denote the complex unit circle.
We want to estimate a signal $x=(x_n)_{n\in V}$, with values $x_n \in \mathbb{S}$, 
%Let us consider a signal $y=(y_n)_{n\in V}$  %$having $N$ elements, for some $N\geq 1$, 
defined on a connected undirected graph $(V,E)$, %with $N$ nodes, for some $N\geq 1$, 
where %$V=\{1,2,\ldots,N\}$  
$V$
is the set of nodes %indexes 
and $E$ is the set of edges, which are sets of two distinct nodes. 
Typically, we are given a noisy signal $y=(y_n)_{n\in V}$ defined on the same graph and the sought signal $x$ is a smoothed, or denoised, version of $y$, which achieves a tradeoff between closeness to $y$ and smoothness, in some sense. Another typical setting is interpolation, or inpainting: $y$ is defined on a subset $U\subset V$ of nodes and we want to estimate its missing samples; that is, $x$ is the smoothest signal defined on $V$ such that $x_n=y_n$, for every $n\in U$.

\subsection{Classical Tikhonov Regularization}  

For %classical 
real-valued signals, Tikhonov-regularized smoothing consists in solving the following convex optimization problem.
%is a simple way to achieve this goal. 
Given $y=(y_n)_{n\in V}$ and nonnegative weights $(w_n)_{n\in V}$ and $(\lambda_{n,n'})_{\{n,n'\}\in E}$, $x=(x_n)_{n\in V}$ is the solution to 
\begin{equation}
%\minimize_{x_n \in \mathbb{S}\,:\, n\in V}
\minimize_{x_n \in \mathbb{R}\,:\, n\in V}
 \,\sum_{n\in V} \frac{w_n}{2} (x_n-y_n)^2 + \! \!\sum_{\{n,n'\}\in E}\! \frac{\lambda_{n,n'}}{2} (x_n-x_{n'})^2.\label{eq1}
\end{equation}
For the interpolating task with $y$ defined only on $U\subset V$, we want to solve, instead:
\begin{equation}
\minimize_{x_n \in \mathbb{R}\,:\, n\in V}  \!\sum_{\{n,n'\}\in E}\! \frac{\lambda_{n,n'}}{2} (x_n-x_{n'})^2\quad\mbox{s.t.}\quad x_n=y_n,\ \forall n\in U.\label{eq2}
\end{equation}
Formally, \eqref{eq2} can be viewed as a particular case of \eqref{eq1} with $w_n=\{+\infty$ if $n\in U$,  0 otherwise$\}$, so that we can focus on the form  \eqref{eq1}, with the weights $w_n$ allowed to be $+\infty$.

We want to formulate an equivalent problem to \eqref{eq1} for signals $x$ and $y$ with values in $\mathbb{S}$. Let us define the argument function $\mathrm{arg}$, which maps $z\in\mathbb{S}$ to $\mathrm{arg}(z)\in(-\pi,\pi]$, such that $z=e^{j\mathrm{arg}(z)}$. A natural adaptation to circle-valued signals is to replace the squared Euclidean distance $(t,t') \in \mathbb{R}^2 \mapsto (t-t')^2$ by the geodesic distance $(z,z') \in \mathbb{S}^2 \mapsto \min ( |\mathrm{arg}(z)-\mathrm{arg}(z')|, 2\pi-|\mathrm{arg}(z)-\mathrm{arg}(z')|)$.  \bl{This yields a  nonconvex and nonsmooth, therefore very difficult, optimization problem to solve. In this work, we consider instead a statistical view of the estimation problem, which leads to a different formulation.}

\subsection{Bayesian View} 

We can notice that \eqref{eq1} corresponds to the maximum-a-posteriori (MAP) estimate of an unknown signal $x^\sharp$ given $y$, which is $x^\sharp$ plus white Gaussian noise, assuming a Gaussian Markov Random Field  prior for $x^\sharp$, with nonzero dependencies between its Gaussian variables along the edges of $V$. That is, $y_n -x^\sharp_n\sim \mathcal{N}(1/w_n)$ and  $x^\sharp_n-x^\sharp_{n'} \sim \mathcal{N}(1/\lambda_{n,n'})$, where $\mathcal{N}(\sigma^2)$ denotes the normal distribution with zero mean and variance $\sigma^2$. Thus, in the circle-valued case, let us consider that $y_n=e^{j\alpha_n}$ is a noisy version of $x^\sharp_n=e^{j\omega^\sharp_n}$, in the sense that $\alpha_n \in \mathbb{R}$ is $\omega^\sharp_n=\mathrm{arg}(x^\sharp_n)\in (-\pi,\pi]$ plus Gaussian noise. Then $\mathrm{arg}(y_n)$, which is the wrapped version in $(-\pi,\pi]$ of $\alpha_n \in \mathbb{R}$, follows the wrapped normal distribution with mean $\omega^\sharp_n$. Since its probability density function (p.d.f.) does not have a closed form, it is common in directional statistics to consider instead, as a close approximation, the von Mises distribution \cite{kha77}. 
That is, we consider that $\mathrm{arg}(y_n)$ is the outcome of a random variable with p.d.f.\ $\propto e^{w_n\cos(\cdot - \omega^\sharp_n)}$. Another argument for the von Mises distribution is that it is the maximum entropy distribution with prescribed `variance' $1/w_n$. Hence, we formulate Tikhonov smoothing for circle-valued signals as the MAP estimate of an unknown Markov Random Field with von Mises dependencies, perturbed by von Mises noise. That is, taking the negative logarithm of the p.d.f., $x_n=e^{j\omega_n}$, where the $\omega_n =\arg(x_n)$ are the solutions to
%$x$ is the solution to 
\begin{align}
\minimize_{\omega_n \in (-\pi,\pi]\,:\, n\in V}
&\,\sum_{n\in V} w_n \big(1-\cos(\omega_n-\mathrm{arg}(y_n))\big)\notag\\
&+ \! \!\sum_{\{n,n'\}\in E}\! \lambda_{n,n'} \big(1-\cos(\omega_n-\omega_{n'})\big) .\label{eq3}
\end{align}
Note that the Taylor series of $1-\cos(t)$ is $t^2/2 + o(t^2)$, so that for small deviations, the problems \eqref{eq3} and \eqref{eq1} behave similarly.

\subsection{Proposed Model}

Now, we can express the problem \eqref{eq3} with respect to the variables $x_n\in \mathbb{S}$, instead of reasoning on their arguments $\omega_n$. The problem becomes:
\begin{align}
%\minimize_{x=(x_n)_{n\in V}} 
\minimize_{x_n \in \mathbb{S}\,:\, n\in V}
&\,\sum_{n\in V} w_n \big(1-\Re(x_n y_n^*)\big) \notag\\
&+ \! \!\sum_{\{n,n'\}\in E}\! \lambda_{n,n'} \big(1- \Re(x_n x_{n'}^*)\big),\label{eq4}
\end{align}
where $\Re$ %$\Re(\cdot)$ 
denotes the real part and $\cdot^*$ denotes the complex conjugation. Note that this problem is nonconvex for two reasons: the variables $x_n$ are constrained to live in the nonconvex circle $\mathbb{S}$ and the product $\Re(x_n x_{n'}^*)$ is nonconvex. The second issue can be resolved by noticing that $1-\Re(x_n x_{n'}^*)=\frac{1}{2}|x_n -x_{n'}|^2$. Indeed $|x_n -x_{n'}|^2=|x_nx_{n'}^*-1|^2=(\Re(x_n x_{n'}^*)-1)^2+\Im(x_n x_{n'}^*)^2=2-2\Re(x_n x_{n'}^*)$, where $\Im$ denotes the imaginary part. Therefore, the problem \eqref{eq4} can be rewritten as:
\begin{equation}
\minimize_{x_n \in \mathbb{S}\,:\, n\in V}\;\sum_{n\in V} \frac{w_n}{2} |x_n-y_n|^2 + \! \!\sum_{\{n,n'\}\in E}\! \frac{\lambda_{n,n'}}{2} |x_n-x_{n'}|^2,\label{eq5}
\end{equation}
which is the natural extension of \eqref{eq1} in the complex plane, and where the objective function to minimize is convex; there remains the nonconvex circle constraint.

Another motivation for our model \eqref{eq4}--\eqref{eq5} is the following: suppose that $y$ is a corrupted version of the unknown circle-valued signal $x^\sharp$ with complex Gaussian noise; that is, independent Gaussian noise with variance $1/w_n$ is added to the real and imaginary parts of each $x^\sharp_n \in \mathbb{S}$. Then, almost surely, the $y_n$ are no longer in $\mathbb{S}$. In that case, the MAP estimate of $x^\sharp$ is exactly the solution to  \eqref{eq5}, where the $y_n$ are now any complex numbers. Moreover,  $\frac{1}{2} |x_n-y_n|^2=\frac{1}{2} (1+|y_n|^2)-\Re(x_n y_n^*) $, so that the problem \eqref{eq5} can be rewritten as:
\begin{align}
%\minimize_{x=(x_n)_{n\in V}} 
\minimize_{x_n \in \mathbb{S}\,:\, n\in V}\; \Psi_\mathrm{orig}(x)=
&\,\sum_{n\in V} w_n \big({\textstyle\frac{1}{2} (1+|y_n|^2)}-\Re(x_n y_n^*)\big) \notag\\
&+ \! \!\sum_{\{n,n'\}\in E}\! \lambda_{n,n'} \big(1- \Re(x_n x_{n'}^*)\big).\label{eq6}
\end{align}
 The problem \eqref{eq6} generalizes \eqref{eq4} to  any complex numbers $y_n$, but since a constant value in the cost function to minimize does not change the solution, \eqref{eq4} and \eqref{eq6} are equivalent.
 
 A natural idea to make the problem \eqref{eq5} convex is to remove the nonconvex circle constraint;
we discuss this approach in Section \ref{sects}. 
Instead, we propose to stick with the formulation \eqref{eq6} and we propose a new convex relaxation of this problem, in Section \ref{secpro}.

\subsection{Limit Cases}\label{seclimc}

Let us look at the two limit cases of \eqref{eq6}, where the data-fit term overwhelms the regularization term, or the other way around. So, let us assume that all $w_n$ are positive and that all $\lambda_{n,n'}$ tend to zero. In that case, the solution is simply $x_n=\{y_n/|y_n|$ if $y_n \neq 0$, any point in $\mathbb{S}$ otherwise$\}$, for every $n\in V$. The other limit case is more interesting: let us assume that all $\lambda_{n,n'}$ tend to $+\infty$. In that case, the regularization term is minimized, which means that the signal is constant: there exists $\mathrm{x}\in\mathbb{S}$ such that $x_n=\mathrm{x}$, for every $n\in V$. This point $\mathrm{x}$ minimizes the data-fit term $\sum_{n\in V} -w_n \Re(\mathrm{x} y_n^*)=-\Re(
\mathrm{x}\sum_{n\in V} w_n y_n^*)=-\Re(\mathrm{x}\mathrm{x}_{\mathrm{av}}^*)$, were $\mathrm{x}_{\mathrm{av}}=\sum_{n\in V} w_n y_n$. The solution, \bl{called the (weighted) circular mean of the points}, is $\mathrm{x}=\{\mathrm{x}_{\mathrm{av}}/|\mathrm{x}_{\mathrm{av}}|$ if $\mathrm{x}_{\mathrm{av}} \neq 0$, any point in $\mathbb{S}$ otherwise$\}$. That is, $\mathrm{x}$ is simply the weighted average of the $y_n$, rescaled to be in $\mathbb{S}$. \bl{Thus, when $w_n\equiv 1$, we have just recovered the well known property that the circular mean of a set of points on the circle is the maximum likelihood estimate of the mean for a von Mises distribution fitting the points.}

\subsection{\bl{Related Work}}\label{sects}

\bl{There is a large literature about optimization on manifolds. For instance, Bergmann and Tenbrinck \cite{ber18} proposed a generic approach for smoothing manifold-valued signals on graphs. Weinmann et al.~\cite{wei14} proposed a method for a large class of functionals including the TV and Tikhonov costs, later extended to the Mumford--Shah functional for piecewise smooth reconstruction \cite{wei16} and to the more general setting of inverse problems \cite{sto21}. The Potts model can be used for the recovery of piecewise-constant signals \cite{wei16}. In general, such approaches are heuristic and have no convergence guarantees; when an algorithm is proved to converge, this is typically to a local solution. Here we focus on estimating the exact global minimizer of the nonconvex problem \eqref{eq6}. For this, we  propose a convex relaxation of the problem, which is tight enough for its solution to coincide, not always but often in practice, with the one of the original problem; when this is the case, this can be certified.}\medskip

\bl{Besides variational regularization, methods based on local averaging can be used. For instance, median filtering for circle-valued data has been proposed \cite{sto18}, which is robust to outliers. 
In the rest of this section, we focus on the} \emph{structure tensor} \cite{zen86,for87,kas87,knu89,big91} a popular tool in image processing to analyze and process the local orientation of a vector field, typically the gradient field of an image. Identifying 2-D vectors with complex numbers, let us first recall that $z\in\mathbb{C}$ and $-z$ have same orientation. Defining the orientation as the angle $\omega=\arg(z)$, we can either restrict $\omega$ to the interval $[0,\pi)$, or define it on the real line modulo $\pi$; that is $\omega\in\mathbb{R}$  is the same orientation as $\omega+\pi$. Therefore, when processing orientations, to avoid cancellation effects that might happen when averaging numbers similar to $z$ with numbers similar to $-z$, it is better to multiply $\omega$ by two, so that it is in $(-\pi,\pi]$, like a phase or angle; that is, one deals with $e^{j2\omega} =z^2$ instead of $z$. Thus, Tikhonov smoothing for orientations consists in solving \eqref{eq4} with the $y_n\in\mathbb{S}$ replaced by $y_n^2$% and $w_n$ replaced by $w_n/4$
; then   the smoothed values  are the $\sqrt{x_n}$, to reverse the doubling operation.

Keeping this squaring effect for orientations in mind,  the structure tensor method, to smooth a 2-D vector field  identified with a complex-valued image $y$, works as follows: each value $y_n\in\mathbb{C}$ is mapped to the $2\times 2$ real matrix of rank 1
\begin{equation}
M_n = \left[ \begin{array}{cc}
\Re(y_n)^2&\Re(y_n)\Im(y_n)\\
\Re(y_n)\Im(y_n)&\Im(y_n)^2
\end{array}\right].\label{eq7}
\end{equation}
Note that $y_n$ and $-y_n$ are mapped to the same matrix, which is consistent with the discussion above: in this context, only the vector orientations matter, not their directions. Then the matrices $M_n$ are spatially averaged by applying a lowpass filter to them, 
elementwise. After filtering, the smoothed matrix $\overline{M}_n$ does not have rank 1, %any more, 
in general. Thus, the smoothed value $x_n$ is obtained by setting $[\Re(x_n)\ \Im(x_n)]^\mathrm{T}$ as the principal eigenvector of $\overline{M}_n$.
To understand 
this process, let us define $\alpha_n=\arg(y_n)$. We have, for every $n\in V$,
\begin{align}
M_n &= |y_n|^2\left[ \begin{array}{cc}
\cos(\alpha_n)^2&\cos(\alpha_n)\sin(\alpha_n)\\
\cos(\alpha_n)\sin(\alpha_n)&\sin(\alpha_n)^2
\end{array}\right]\\
&= \frac{|y_n|^2}{2}\left[ \begin{array}{cc}
1&0\\
0&1
\end{array}\right]+\frac{|y_n|^2}{2}\left[ \begin{array}{cc}
\cos(2\alpha_n)&\sin(2\alpha_n)\\
\sin(2\alpha_n)&-\cos(2\alpha_n)
\end{array}\right].\notag
\end{align}
After filtering, $\overline{M}_n$ is symmetric and can be written as 
%and positive semidefinite, of the form
\begin{align}
\overline{M}_n &= \left[ \begin{array}{cc}
m_{n,1}&m_{n,3}\\
m_{n,3}&m_{n,2}
\end{array}\right]\\
&= d_{n}\left[ \begin{array}{cc}
1&0\\
0&1
\end{array}\right]+a_n\left[ \begin{array}{cc}
\cos(2\omega_n)&\sin(2\omega_n)\\
\sin(2\omega_n)&-\cos(2\omega_n)
\end{array}\right]\notag\\
&=\left[\!\!\!\begin{array}{cc}
2a_n\cos(\omega_n)^2\!+\!d_{n}\!-\!a_n\!&\!2a_n\cos(\omega_n)\sin(\omega_n)\\
2a_n\cos(\omega_n)\sin(\omega_n)\!&\!2a_n\sin(\omega_n)^2\!+\!d_{n}\!-\!a_n
\end{array}\!\!\!\right]\notag\\
&=\left[\!\!\begin{array}{cc}
\cos(\omega_n)&\sin(\omega_n)\\
\sin(\omega_n)&-\cos(\omega_n)
\end{array}\!\!\right]\left[\!\!\begin{array}{cc}
d_{n}+a_n&0\\0&d_{n}-a_n
\end{array}\!\!\right]\notag\\
&\quad\ \times\left[\!\!\begin{array}{cc}
\cos(\omega_n)&\sin(\omega_n)\\
\sin(\omega_n)&-\cos(\omega_n)
\end{array}\!\!\right],\notag
\end{align}
with $d_{n}=(m_{n,1}+m_{n,2})/2$, $a_n=\big((m_{n,1}-m_{n,2})^2/4+m_{n,3}^2\big)^{1/2}=|x_n|/2$, and $\omega_n=\arg(x_n)/2$, where we set 
\begin{equation}
x_n=(m_{n,1}-m_{n,2})+2j m_{n,3}=2a_n e^{j2\omega_n}.
\end{equation}
%$\omega_n\in [0,\pi)$,  
Thus, the two eigenvalues of $\overline{M}_n$, in decreasing order, are $d_{n}+a_n$ and $d_{n}-a_n$ and its principal eigenvector is $[\cos(\omega_n)\ \sin(\omega_n)]^\mathrm{T}$. $a_n$ is a confidence indicator: if $a_n=0$, there is no preferred direction locally, whereas if it is large, the direction $\omega_n$ is dominant. Thus, we \bl{do not} need the matrix formalism: reasoning on the complex numbers $y_n$ and $x_n$ is equivalent and easier.  
Indeed, since $y_n^2 = |y_n|^2 e^{j2\alpha_n}=(\Re(y_n)^2-\Im(y_n)^2)+2j\Re(y_n)\Im(y_n)$, $x_n$ is simply the result of spatial averaging applied to the $y_n^2$.

Tikhonov regularization amounts to lowpass filtering: with $w_n\equiv 1$ and $\lambda_{n,n'}\equiv \lambda$, the solution $x$ to \eqref{eq5} without the circle constraint is simply the result of a convolution applied to $y$, with inverse frequency response one plus $\lambda$ times the graph Laplacian. Therefore, the structure tensor method is essentially solving the Tikhonov problem \eqref{eq5}, with the variables searched in $\mathbb{C}$ instead of $\mathbb{S}$ and with the $y_n$ replaced by $y_n^2$, this second change being specific to the setting of orientations. Note that if $|y_n|=1$ for all $n$, squaring the $y_n$ in \eqref{eq5}, with the circle constraint, still corresponds to a MAP estimate with a scaled von Mises prior. However, if the amplitudes $|y_n|$ are arbitrary, squaring the $y_n$ also squares their amplitudes, so that the squared amplitudes are averaged by the regularization process; there seems to be no obvious Bayesian interpretation of \eqref{eq5}, with or without the circle constraint, in that case. Thus, it is better to divide $M_n$ by $|y_n|$ in \eqref{eq7} before spatial averaging, which is the way the structure tensor is defined by Knutsson \cite{knu89}.

Finally, let us remark that if $|y_n|=1$, when solving \eqref{eq5} without the circle constraint, 
the $x_n$ remain in the convex hull of the $y_n$, which is contained in the 
complex disk $\mathbb{D}=\{z\in\mathbb{C}\ :\ |z|\leq 1\}$, the convex hull of $\mathbb{S}$. 
 Thus, there is no need to enforce the constraint that the $x_n$ belong to $\mathbb{D}$, since it is automatically satisfied. In the sequel, we will refer to the following process, to find an approximate solution of \eqref{eq5}, as the \emph{baseline method}: \eqref{eq5} is solved without the circle constraint (which amounts to solving a linear system) and the $x_n$ are rescaled as $x_n/|x_n|$ afterwards, to make them lie in $\mathbb{S}$.

\section{Fourier Lifting: Convex relaxation using moments of measures}\label{secpro}

\textbf{The method of moments} -- There is a general recipe to reformulate, or \emph{lift}, a nonconvex problem as a convex one: the minimization of a function $f$ over $\mathbb{S}$ is equivalent to minimizing $\int_{\mathbb{S}} f(z)\mathrm{d}\mu(z)$ over $\mu$ in the set of probability measures (i.e. positive Borel measure with mass 1) on $\mathbb{S}$, assuming that $\arg\min(f)$, the set of minimizers of $f$, is nonempty. 
This latter problem is convex, since it consists in minimizing a linear functional over a convex set. A minimizing measure $\mu^\star$ will be concentrated over $\arg\min(f)$; in particular, if the minimizer $z^\star$ of $f$ is unique, $\mu^\star=\delta_{z^\star}$, the Dirac measure at  $z^\star$. This principled approach has a major downside, yet: the set of probability measures is infinite-dimensional, which prevents its numerical implementation in  general. However, there is a case where the method can be implemented exactly: if the measure can be parameterized and recovered from a finite number $M$ of its moments $\hat{\mu}_m=\int_{\mathbb{S}} \phi_m(z) \mathrm{d}\mu(z)$, $m=1,\ldots,M$, for some basis functions  $\phi_m$, and if $f=\sum_{m=1}^M a_m \phi_m$ is a linear combination of the $\phi_m$, then $\int_{\mathbb{S}} f(z)\mathrm{d}\mu(z)= \sum_{m=1}^M a_m \hat{\mu}_m$, so that the problem becomes convex and finite-dimensional, in terms of the moments $\hat{\mu}_m$: we want to minimize the linear term $\sum_{m=1}^M a_m c_m$ with respect to the coefficients $c_m$, under the constraint that $c_m=\hat{\mu}_m$ for every $m$, for some probability measure $\mu$. This approach is called \emph{the method of moments} \cite{mez03}. 
In this work, we use trigonometric moments, or Fourier coefficients. The characterization of a measure on the circle from a subset of its Fourier coefficients has a long history, rooted in Carath\'eodory's work~\cite{car11}; \bl{related theorems are often called Bochner's  theorems}. In short, the constraint that the $c_m=\hat{\mu}_m$ for some positive measure $\mu$ on the circle 
 is satisfied if the Toeplitz matrix formed by the $c_m$ is positive semidefinite \cite{cur91,con20}.\medskip

\textbf{Convex relaxation of optimization over graphs using measures} --  Our problem \eqref{eq6} features nonconvex pairwise costs $(x_n,x_{n'}) \mapsto  \Re(x_n x_{n'}^*)$. The minimization of a function $g(z,z')$ is equivalent to minimizing $\int_{\mathbb{S}^2} g(z,z')\mathrm{d}\nu(z,z')$ over $\nu$ in the set of probability measures on $\mathbb{S}^2$. Hence, given  an optimization problem  over the graph $(V,E)$ with unary potential costs $f_n$ at the nodes and symmetric pairwise interaction costs $g_{n,n'}$ at the edges, all bounded from below and lower semicontinuous:
\begin{equation}
\minimize_{x_n \in \mathbb{S}\;:\; n\in V} \;\sum_{n\in V} f_n (x_n) + \! \!\sum_{\{n,n'\}\in E}\! g_{n,n'}(x_n, x_{n'}),\label{eq10}
\end{equation}
we propose the following lifting technique: we introduce a probability measure $\mu_n$ on $\mathbb{S}$ for each node $n\in V$, as well as a probability measure $\nu_{n,n'}$ on $\mathbb{S}^2$ for each edge $\{n,n'\}\in E$, and we consider the lifted convex problem:
\begin{align}
\minimize_{(\mu_n),(\nu_{n,n'})} &\;\sum_{n\in V}\int_{\mathbb{S}} f_{n}(z)\mathrm{d}\mu_{n}(z)\notag\\
&+ \! \!\sum_{\{n,n'\}\in E}\
\int_{\mathbb{S}^2} g_{n,n'}(z,z')\mathrm{d}\nu_{n,n'}(z,z')\label{eq11}\\
&\mbox{s.t. the two marginals of } \nu_{n,n'} \mbox{ are } \mu_{n} \mbox{ and } \mu_{n'}.\notag
\end{align}
It is important to note that, in general, this convex relaxation is not tight: even if the solution $(x^\star_n)_{n\in V}$ to \eqref{eq10} is unique, 
it is not guaranteed that the solution to \eqref{eq11} corresponds to $\mu^\star_n=\delta_{x_n}$: it might be that measures which are not Diracs achieve a lower value of the objective function. 

Interestingly, there is a strong connection with the theory of optimal transport \cite{vil03}: the function 
\begin{align}
(\mu_n,\mu_{n'}) \mapsto &\min_{\nu_{n,n'}} \int_{\mathbb{S}^2} g_{n,n'}(z,z')\mathrm{d}\nu_{n,n'}(z,z')\label{eq12}\\
&\mbox{s.t. the two marginals of } \nu_{n,n'} \mbox{ are } \mu_{n} \mbox{ and } \mu_{n'},\notag
\end{align}
is the Monge--Kantorovich optimal transport cost between the probability measures $\mu_{n}$ and  $\mu_{n'}$, interpreting $g_{n,n'}(z,z')$ as the cost of moving one unit of mass from the point $z$ to the point $z'$; the minimizing measure $\nu^\star_{n,n'}$ in \eqref{eq12}, which exists, is called the optimal coupling measure.  We refer to \cite{rab11} for more details on optimal transport on the circle $\mathbb{S}$.

When the measures are restricted to live on a finite set of labels, instead of a continuous set like $\mathbb{S}$, the relaxation \eqref{eq11}, which is a linear program, is well known in statistics, in the fields of graphical models, discrete inference and labeling, where it is called the local polytope relaxation \cite{wer07,wai08,kap15}. In the present work, we \bl{do not} want to discretize the circle $\mathbb{S}$. We will instead parameterize the measures by a finite number of their Fourier coefficients, like in the method of moments; hence, we name our convex relaxation approach as \emph{Fourier lifting}.\medskip

\textbf{Fourier lifting} --  
A probability measure $\nu$ defined on $\mathbb{S}^2$ has Fourier coefficients 
$\hat{\nu}_{m,m'}=\int\!\!\int_{\mathbb{S}^2}z^{-m} {z'}^{-m'}\mathrm{d}\nu(z,z')$
%$\hat{\nu}_{m,n}=\int_{[-\pi,\pi)^2} e^{-j(m\omega_1+n\omega_2)}\nu(e^{j\omega_1},e^{j\omega_2}) d(\omega_1,\omega_2)$, 
for every
$(m,m')\in\mathbb{Z}^2$. 
We have $\hat{\nu}_{-m,-m'}=\hat{\nu}_{m,m'}^*$%, for every $(m,n)$, 
 and $\hat{\nu}_{0,0}=1$.

$\nu$ has two marginals $\mu=\int_\mathbb{S}\mathrm{d}\nu(\cdot,z)$ and $\mu'=\int_\mathbb{S}\mathrm{d}\nu(z,\cdot)$, which are probability measures on $\mathbb{S}$. They have Fourier coefficients 
$\hat{\mu}_{m}=\int_{(-\pi,\pi]} e^{-jm\omega}\mu(e^{j\omega})d\omega=\hat{\nu}_{m,0}$ and  $\hat{\mu}'_{m'}=\int_{(-\pi,\pi]} e^{-jm'\omega}\mu'(e^{j\omega})d\omega=\hat{\nu}_{0,m'}$ respectively, for every 
$(m,m')\in\mathbb{Z}^2$.

If $\nu$ is a 2-D Dirac in $(e^{j\omega_1},e^{j\omega_2})\in\mathbb{S}^2$, $c_{m,m'}=e^{-j(m\omega_1+n\omega_2)}=c_{m,0}c_{0,m'}$, so that the matrix of moments has rank 1 and $|c_{m,m'}|= 1$ for every $(m,m')$.

In this work, we will only parameterize $\nu$ using $c_{0,0}=1$, $c_{1,0}$, $c_{0,1}$, $c_{-1,1}$. Indeed, if $\nu$ is a Dirac in $(e^{j\omega_1},e^{j\omega_2})$, $c_{1,0}=e^{-j\omega_1}$, $c_{0,1}=e^{-j\omega_2}$, and 
\begin{align}
\cos(\omega_1-\omega_2) &= \int \cos(\omega-\omega') \nu(e^{j\omega},e^{j\omega'}) d(\omega,\omega') \\
&= \frac{1}{2}(c_{-1,1}+c_{1,-1})=\Re(c_{-1,1}).
\end{align}

Let us map these coefficients in the $3\times 3$ Hermitian matrix
\begin{equation}
P=\left[\begin{array}{ccc}
1 &c_{1,0}&c_{0,1} \\
c_{1,0}^* &1&c_{-1,1}\\
c_{0,1}^*&c_{-1,1}^*&1
\end{array}\right]
\end{equation}

$P$ is positive semidefinite, which we denote by $P\succcurlyeq 0$, if and only if all its principal minors are nonnegative; that is,
\begin{equation}
\left|\begin{array}{ccc}
1 &c_{1,0}&c_{0,1} \\
c_{1,0}^* &1&c_{-1,1}\\
c_{0,1}^*&c_{-1,1}^*&1
\end{array}\right|\geq 0,\left|\begin{array}{cc}
1 &c_{1,0} \\
c_{1,0}^* &1
\end{array}\right|\geq 0,
\end{equation}
\begin{equation}\left|\begin{array}{cc}
1 &c_{0,1} \\
c_{0,1}^* &1
\end{array}\right|\geq 0,
\left|\begin{array}{cc}
1&c_{-1,1}\\
c_{-1,1}^*&1
\end{array}\right|\geq 0.
\end{equation}
Equivalently, $|c_{1,0}|\leq 1$, $|c_{0,1}|\leq 1$, $|c_{-1,1}|\leq 1$, $1+2\Re(c_{1,0}c_{-1,1}c_{0,1}^*)-|c_{1,0}|^2-|c_{0,1}|^2-|c_{-1,1}|^2\geq 0$.

If $\nu$ is a Dirac, $P\succcurlyeq 0$. Thus, for every probability measure $\nu$, by convexity of the positive semidefinite cone, $P\succcurlyeq 0$. Moreover, $P\succcurlyeq 0$ has rank 1 if and only if  $\nu$ is a Dirac; that is,
\begin{equation}
P=\left[\begin{array}{c}
1  \\
c_{1,0}^* \\
c_{0,1}^*
\end{array}\right]\left[\begin{array}{ccc}
1 &c_{1,0}&c_{0,1}
\end{array}\right],
\end{equation}
with $|c_{1,0}|=1$, $|c_{0,1}|=1$.\\

Hence, the proposed convex relaxation of \eqref{eq6} is:
\begin{align}
&\minimize_{\small\substack{x_n \in \mathbb{C}\;:\; n\in V\\ r_{n,n'}\in \mathbb{C}\;:\; (n,n')\in E}} \; \Psi_\mathrm{conv}(x,r)=\sum_{n\in V}w_n \big({\textstyle\frac{1}{2} (1+|y_n|^2)}-\Re(x_n y_n^*)\big)\notag\\
&\qquad+{}\!\!\!\!\!\! \sum_{\{n,n'\}\in E}\!\!\!\lambda_{n,n'} \big(1- \Re(r_{n,n'})\big) \notag\\ 
&\qquad\mbox{s.t.}\ \left[\begin{array}{ccc}
1 &x_n^*&x_{n'}^* \\
x_n &1&r_{n,n'}\\
x_{n'}&r_{n,n'}^*&1
\end{array}\right]\succcurlyeq 0,\ \forall (n,n')\in E.\label{eq19}
\end{align}

\begin{figure*}[t]
\centering
\hspace*{-2.4cm}\includegraphics[scale=0.75]{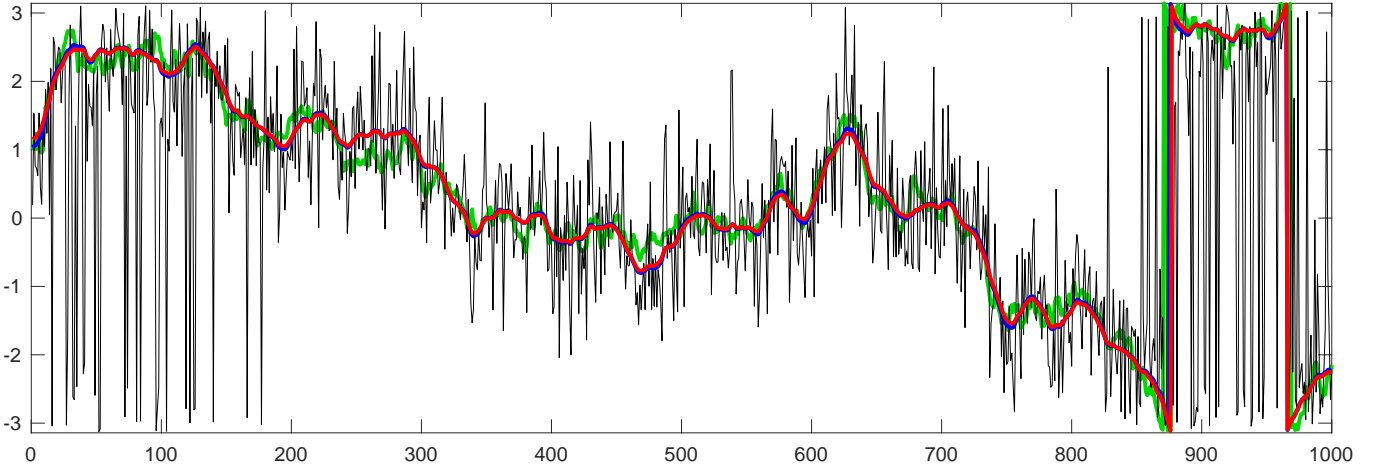}
\caption{Denoising of a circle-valued signal. In green, the ground-truth signal, in black, the noisy signal, in blue, the signal denoised with the baseline method and in red, the signal $(x^\star_n)_{n\in V}$ denoised with the proposed approach, which turns out to be the exact solution to \eqref{eq4}. All signals have their values in $\mathbb{S}$, whose argument in $(-\pi,\pi]$ is displayed.}
\end{figure*}

A Hermitian positive semidefinite matrix with ones on its diagonal is a correlation matrix, the set of which is sometimes called an elliptope \cite{chr79}. An elliptope is convex and compact. It is known that a linear function attains its minimum over a compact convex set at a point of its boundary. Since the rank-1 matrices associated to Diracs are extreme points of the elliptope,  linear minimization over an elliptope is likely to yield a rank-1 matrix. The convex \bl{optimization} problem \eqref{eq19} can be viewed as  linear minimization over a product of elliptopes with linear equality constraints; 
there is no guarantee that all matrices in \eqref{eq19} will be of rank 1 at a solution, but this is what we hope for.

Let us call $(x^\star_n)_{n\in V}$ and $(r^\star_{n,n'})_{(n,n')\in E}$ a solution obtained by solving \eqref{eq19} and $\Psi_\mathrm{conv}^\star$ the corresponding minimal objective value. We also denote by $\Psi_\mathrm{orig}^\star$ the minimal objective value of the original nonconvex problem \eqref{eq6}. 
If all matrices appearing in \eqref{eq19} are of rank 1, or equivalently $|x^\star_n|=1$ for every $n\in $ and $r^\star_{n,n'} = x^{\star}_n x_{n'}^{\star *} $ for every $(n,n')\in E$,  $\Psi_\mathrm{conv}^\star= \Psi_\mathrm{orig}^\star$ and we have obtained an exact solution of the original problem \eqref{eq6}. Otherwise, we rescale the $x^\star_n$ as $x^\star_n/|x^\star_n|$ to project them on $\mathbb{S}$, and we now have an approximate solution $x$ to \eqref{eq6}. Let us denote by $\Psi_\mathrm{approx}=\Psi_\mathrm{orig}(x)$ the objective value evaluated at this $x$. We have $\Psi_\mathrm{conv}^\star\leq \Psi_\mathrm{orig}^\star\leq \Psi_\mathrm{approx}$, so that we can use $(\Psi_\mathrm{approx}-\Psi_\mathrm{conv}^\star)/\Psi_\mathrm{conv}^\star$ as a measure of relative suboptimality of the convex relaxation with respect to the original problem. \bl{We conjecture that the proposed relaxation is tight and yields the exact solution whenever the graph $(V,E)$ has no cycle, as is the case for a 1-D chain.}

\section{Proposed Algorithm}\label{sec4}

We endow $\mathbb{C}$ with the inner product %$(z,z')\in\mathbb{C}^2\mapsto 
$\langle z,z'\rangle = \Re(zz'^*)$, to form a real Hilbert space. Then the problem \eqref{eq19} has the form 
\begin{align}
&\minimize_{s\in \mathbb{C}^d} \; \Psi_\mathrm{conv}(s)=\langle s,e\rangle + f(Ls),\label{eq21}
\end{align}
where the variable $s$ is the concatenation of all $x_n$ and $r_{n,n'}$, the dimension $d$ is the total number of nodes and edges, $e\in \mathbb{C}^d$ is the concatenation of all constants $-w_n y_n$ and $-\lambda_{n,n'}$, the linear operator $L$ maps $e$ to the concatenation of matrices 
\begin{equation}
(Le)_{n,n'}=\left[\begin{array}{ccc}
0 &x_n^*&x_{n'}^* \\
x_n &0&r_{n,n'}\\
x_{n'}&r_{n,n'}^*&0
\end{array}\right]
\end{equation}
for all $(n,n')\in E$, $f: (Q_{n,n'})_{(n,n')\in E}\mapsto \sum_{(n,n')\in E}  \{0$ if $Q_{n,n'}+\mathrm{Id}\succcurlyeq 0$, $+\infty$ otherwise$\}$, and $\mathrm{Id}$ is the $3\times 3$ identity matrix. We endow the set of $3\times 3$ Hermitian matrices with the Frobenius inner product %$(Q,Q')\mapsto 
$\langle Q,Q'\rangle=\mathrm{tr}(QQ')$, where $\mathrm{tr}$ denotes the trace.

To solve the problem \eqref{eq21}, a well suited algorithm is the Proximal Method of Multipliers \cite{roc762,con192}, which, initialized with some variables $U^{(0)}\in (\mathbb{C}^{3\times 3})^{|E|}$ and $s^{(0)}\in\mathbb{C}^d$, consists in the iteration: for $i=0,1,\ldots$
\begin{align*}
%&\,\mbox{for }i=0,1,\ldots\notag\\[-1mm]
    &\left\lfloor
    \begin{array}{l}
    a^{(i)}=L^* U^{(i)}+e\\
    s^{(i+1)}=s^{(i)}- \tau a^{(i)}\\
    U^{(i+1)}=\mathrm{prox}_{\sigma f^*}\big(U^{(i)}+\sigma L (s^{(i+1)}-\tau a^{(i)} ) \big)\\
\end{array}\right.\end{align*}
where $L^*$ denotes the adjoint operator of $L$, $f^*$ denotes the convex conjugate of $f$ \cite{bau17}, $\tau>0$ is a parameter, and we set $\sigma=1/(\|L\|^2\tau)$, where the squared operator norm $\|L\|^2$ is twice the maximum number of edges per node.
 With this choice, the variable $s^{(i)}$ in the algorithm converges to a solution $s^\star$ of \eqref{eq21} \cite[Theorem 4.3]{con192}. In the algorithm, the proximity operator $\mathrm{prox}_{\sigma f^*}$ maps each matrix $Q_{n,n'}$, for $(n,n')\in E$, to the projection of $Q_{n,n'}+\sigma\mathrm{Id}$ onto the cone of Hermitian negative semidefinite matrices, minus $\sigma\mathrm{Id}$; this is achieved by computing the eigendecomposition and setting the positive eigenvalues to zero.

\begin{figure}[t]
\centering
\vspace*{-2mm}\begin{tabular}{cc}
\hspace*{-8mm}\includegraphics[scale=0.445]{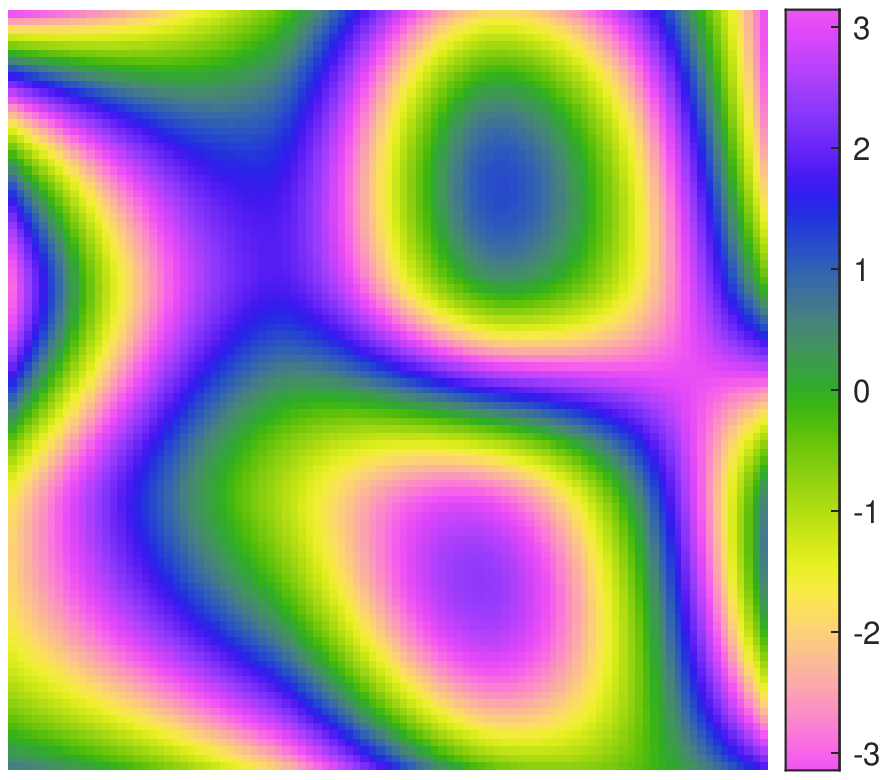}&\hspace*{-5mm}\includegraphics[scale=0.445]{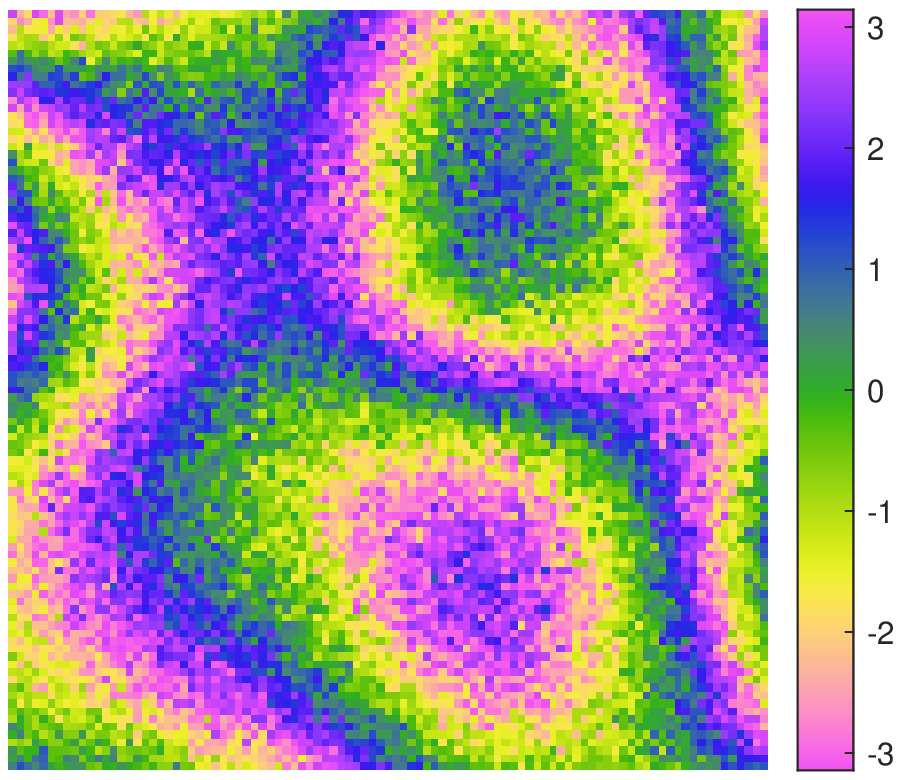}\\[-3mm]
(a)&(b)\\
\hspace*{-8mm}\includegraphics[scale=0.445]{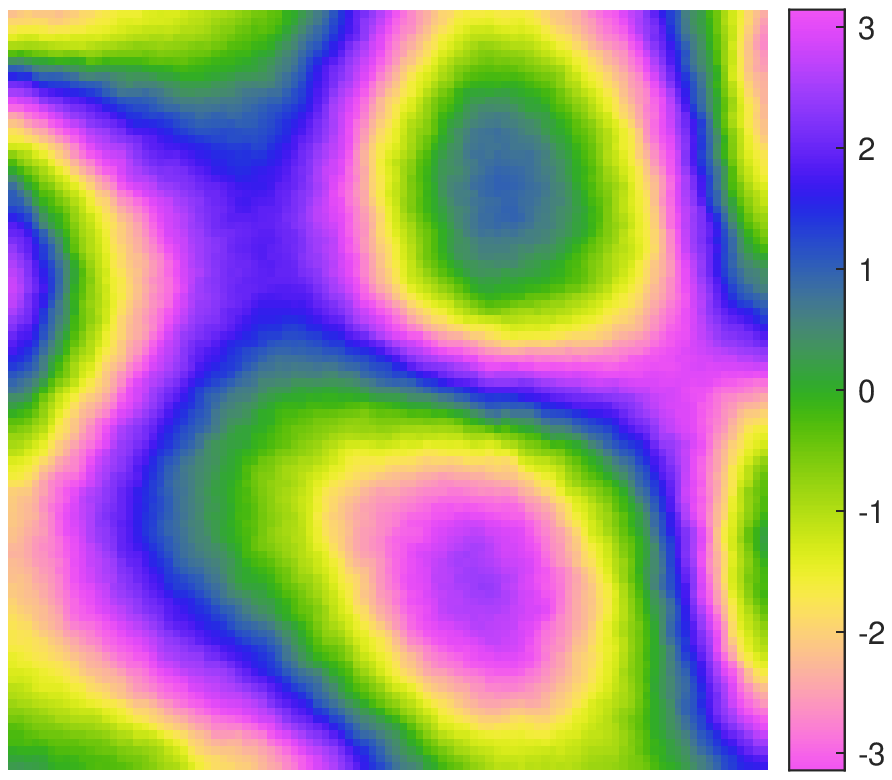}&\hspace*{-5mm}\includegraphics[scale=0.445]{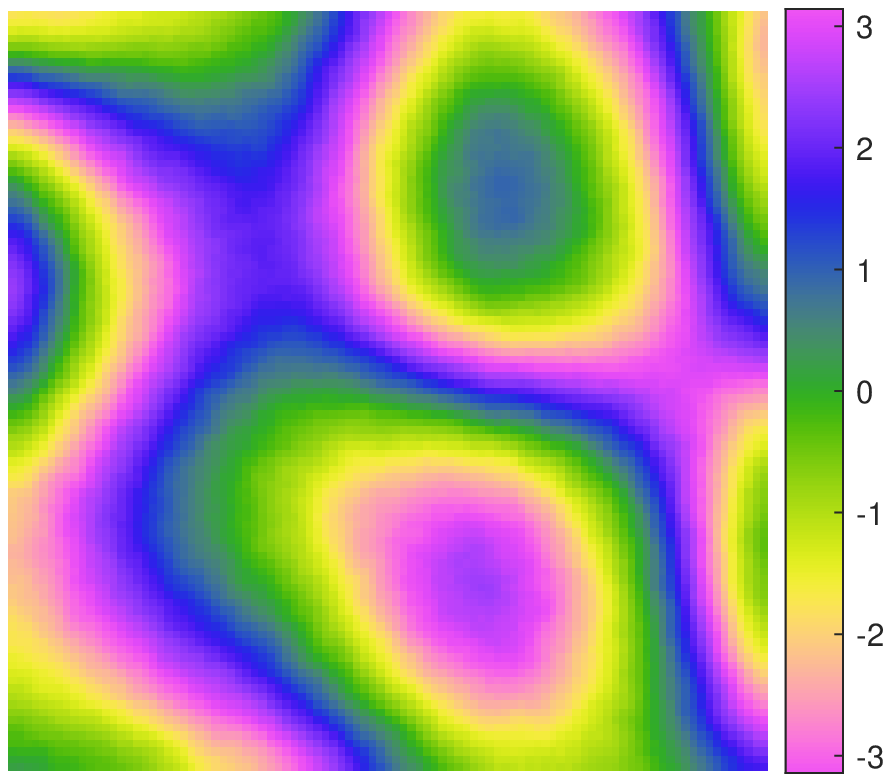}\\[-3mm]
(c)&(d)\\
\hspace*{-8mm}\includegraphics[scale=0.445]{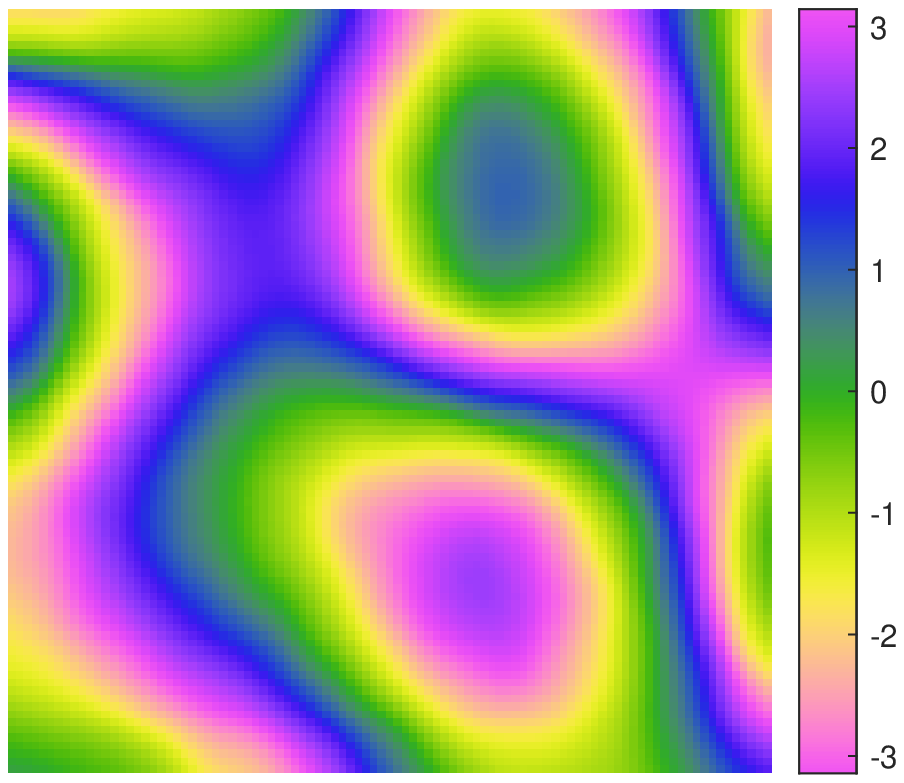}\\[-3mm]
(e)
%$\Psi_\mathrm{orig}(x)=2534$&\!(d) $\Psi_\mathrm{orig}^\star=\Psi_\mathrm{conv}^\star=2479$
\end{tabular}
\caption{Denoising of a circle-valued image: (a) the ground-truth image, (b) the noisy image, (c) the image denoised with the baseline method, (d) the image denoised with the proposed method, which turns out to be the exact solution to \eqref{eq4}, (e) the image denoised by circular mean filtering with Gaussian weights. All images have their values in $\mathbb{S}$, whose argument in $(-\pi,\pi]$ is displayed, using the C2 cyclic colormap designed by Peter Kovesi \cite{kov15}.}
\end{figure}

\begin{figure}[t]
\centering
\vspace*{-2mm}\begin{tabular}{cc}
\hspace*{-7mm}\includegraphics[scale=0.435]{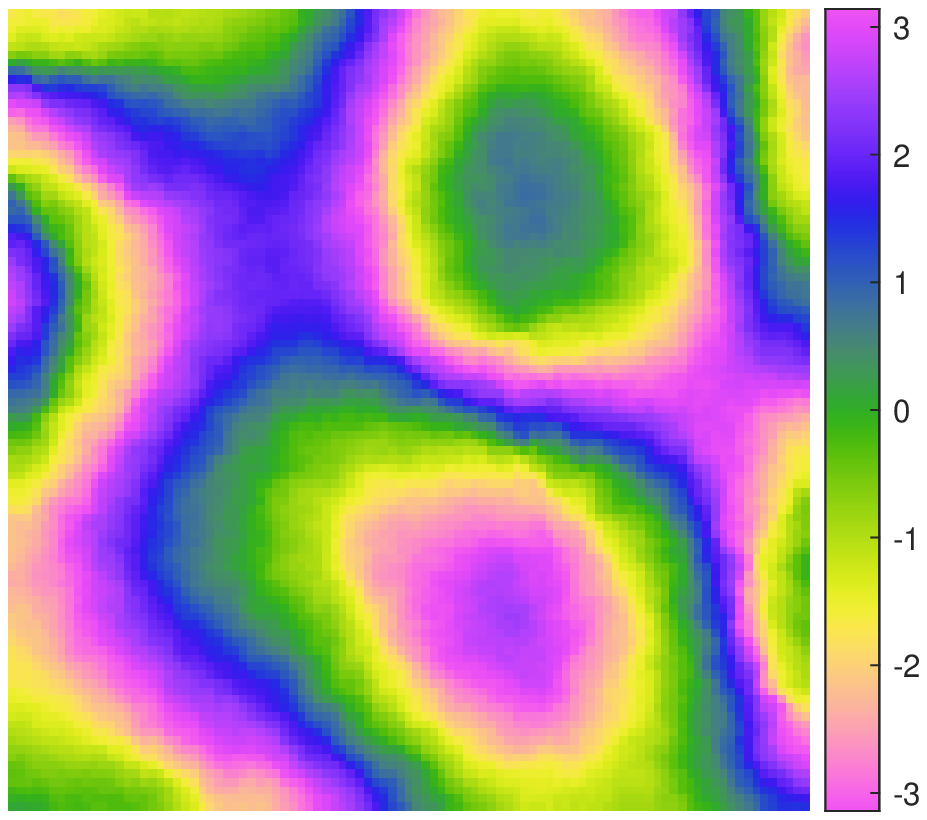}&\hspace*{-4.5mm}\includegraphics[scale=0.435]{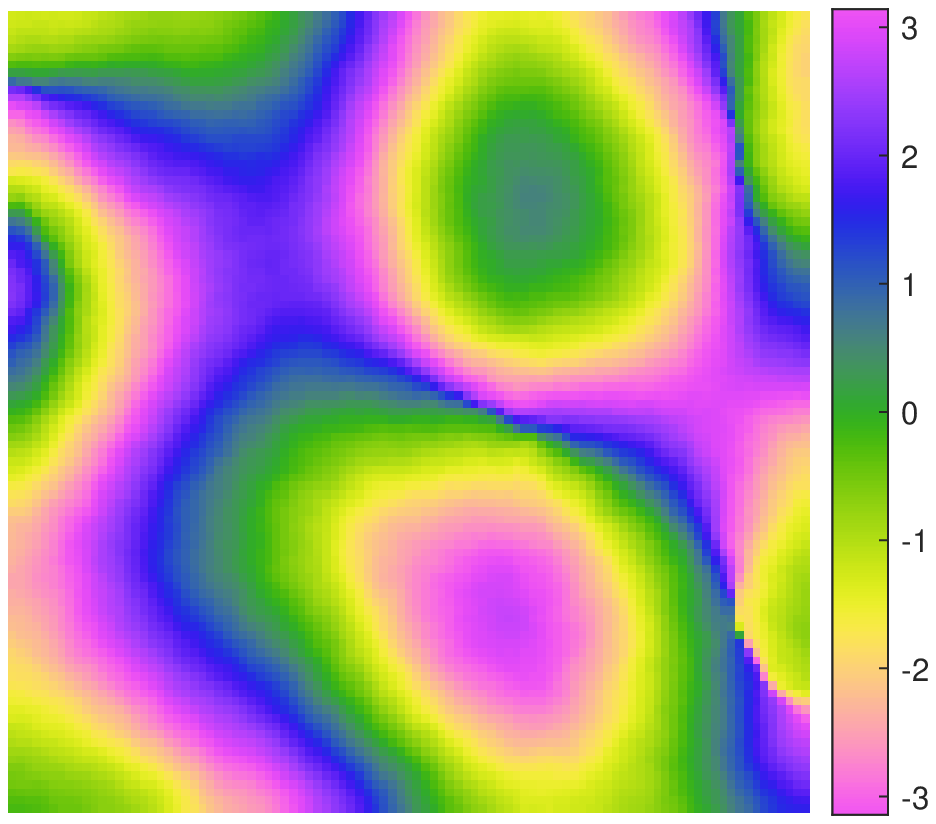}\\[-3mm]
(a)&(b)
\end{tabular}
\caption{Denoising of a circle-valued image, like in Fig.~2, but with the noise level and regularization parameter $\lambda$ twice higher. (a) the image denoised with the baseline method, (b) the image denoised with the proposed method, which is only an approximate solution to  \eqref{eq4}.}
\end{figure}

\section{Experiments}\label{sec5}

For the following experiments, MATLAB code implementing the algorithms and generating the images in the figures is available on the author's webpage. \bl{The code was run in MATLAB R2022a on a Apple Macbook Pro 2019 laptop. }

\subsection{Denoising of a 1-D signal}

In a first experiment, we denoise a 1-D signal of size $N=1000$. That is, $V=\{1,\ldots,N\}$ and $E=\{(1,2),\ldots,(N-1,N)\}$. The gound-truth signal $(e^{j\omega^\sharp_n})_{n\in V}$ is generated using $\omega^\sharp_1=1$ and i.i.d random increments $\omega^\sharp_{n+1}-\omega^\sharp_{n}$ following the Gaussian law of standard deviation $0.1$. Then the noisy signal $y$ is formed by adding to the $\omega_n^\sharp$ white Gaussian noise of 
standard deviation $\sqrt{\lambda}/10$, where $\lambda=50$. These two signals are shown in Fig.~1 in green and black, respectively. We denoise $y$ using the baseline method described at the end of Section~\ref{sec2}, with $\omega_n\equiv 1$ and $\lambda_{n,n'}\equiv \lambda$; the denoised signal is shown in blue in Fig.~1. The corresponding cost value  in \eqref{eq6}, or equivalently in \eqref{eq3}, \eqref{eq4}, or \eqref{eq5}, is $\Psi_\mathrm{orig}(x)\approx 227$. Then we apply the proposed approach by solving \eqref{eq19}, with $\tau=0.1$  in the algorithm, which converges to machine precision in about $300$ iterations. 
 The obtained denoised signal, shown in red in Fig.~1, satisfies $|x^\star_n|=1$ and $r^\star_{n,n+1} = x^{\star}_n x_{n+1}^{\star *} $ for every $n$, so that it is the exact solution to the original problem \eqref{eq6}. The corresponding optimal cost value is $\Psi_\mathrm{conv}^\star= \Psi_\mathrm{orig}^\star\approx 226$. The quantitative and qualitative difference between the results of the baseline and proposed methods is small in this example, but it is satisfying to be able to solve the nonconvex problem of Tikhonov smoothing exactly. 

\subsection{Denoising of a 2-D image}

In a second experiment, we denoise a 2-D image: the phase of a smooth ground-truth image of size $97\times 97$ is generated by cubic interpolation from a random $4\times 4$ image and white Gaussian noise of 
standard deviation 0.5 is added to the phases to obtain a noisy version $y$; they are shown in Fig.~2 (a) and (b), respectively. The graph is the classical square grid: there is a node at each pixel and the edges connect all pairs of horizontally or vertically adjacent pixels. We denoise $y$ using the baseline method described at the end of Section~\ref{sec2}, with $\omega_n\equiv 1$ and $\lambda_{n,n'}\equiv \lambda=5$ \bl{(with $400$ iterations, computation time 0.09s)}; the denoised image is shown in Fig.~2 (c). The corresponding cost value  in \eqref{eq6} is $\Psi_\mathrm{orig}(x)\approx 2534$. Then we apply the proposed approach by solving \eqref{eq19}, with $\tau=0.1$  in the algorithm, which converges to machine precision in about $400$ iterations \bl{(computation time 73s). }
 The obtained denoised image, shown in Fig.~2 (d), satisfies $|x^\star_n|=1$ and $r^\star_{n,n'} = x^{\star}_n x_{n'}^{\star *} $ for every $(n,n')$, so that it is the exact solution to the original problem \eqref{eq6}. The corresponding optimal cost value is $\Psi_\mathrm{conv}^\star= \Psi_\mathrm{orig}^\star\approx 2479$. By comparing the images in Fig.~2 (c) and (d), we can see that the image with the proposed method is a bit more regular, with less jagged  level lines. \bl{We also show in Fig.~2 (e) the image obtained by replacing each pixel value by the weighted circular mean of its neighbors, with Gaussian weights, as described in Section~\ref{seclimc} (computation time 0.003s); that is, we simply apply to $y$ a convolution with a Gaussian filter (of standard deviation 3 pixels) and we rescale each value $x_n$ to project it on the circle. The image is smooth and visually pleasant but its cost value is $\Psi_\mathrm{orig}(x)\approx 2542$, similar to the one of the baseline method.}
 
 \bl{Our current implementation of the proposed algorithm is slow but it calls the eigendecomposition of every matrix to project it on the cone of positive semidefinite matrices. A careful implementation with a routine dedicated to this projection for $3\times 3$ Hermitian matrices would reduce the computation time significantly.}
 
 It is not always the case that the proposed convex relaxation is tight and yields the solution to the original problem. For instance, if we keep the same experiment but with a noise standard deviation of 1 and  $\lambda=10$, the solution of \eqref{eq19} does not satisfy $|x^\star_n|=1$ for all $n$ any more. We have $\Psi_\mathrm{approx}=6797>\Psi_\mathrm{conv}^\star\approx 5953$. The baseline method yields an image $x$ with a cost of $\Psi_\mathrm{orig}(x)\approx 6749$. Both images are shown in Fig.~3. In particular, we can see that the proposed method has introduced an incorrect junction at the bottom right of the image in Fig.~3 (b). Elsewhere, it is more regular and probably closer to the exact solution of the problem \eqref{eq6} than with the baseline method. Thus, the proposed approach is best suited when the noise level is not too high.

\subsection{Interpolation of a 1-D signal}

\begin{figure}[t]
\centering
\hspace*{-6mm}\includegraphics[scale=0.75]{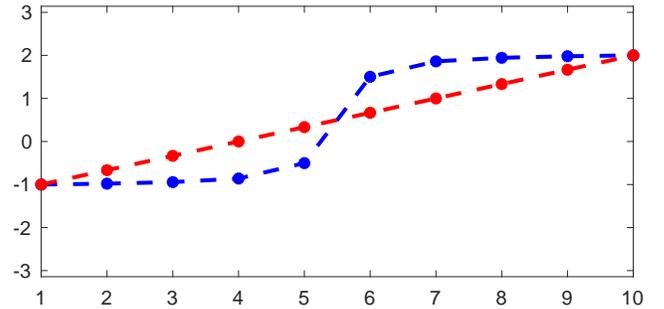}
\caption{Interpolation at $n=2,\ldots,9$ between the fixed boundary values $y_1=e^{-j}$ and $y_{10}=e^{2j}$. In blue, with the baseline method; in red, with the proposed method, which yields the expected solution; that is, a straight line.}
\end{figure}

We now consider interpolation at the intermediate indexes $n=2,\ldots, 9$ of the 1-D signal $y$ defined at $n=1$ and $n=10$ by $y_1=e^{-j}$ and $y_{10}=e^{2j}$. The problem we would like to solve is
\begin{equation}
\minimize_{x_n \in \mathbb{S}\,:\, n=1,\ldots,10}\;\sum_{n=1}^9{\textstyle \frac{1}{2}} |x_{n+1}-x_{n}|^2\ \ \mbox{s.t.}\ \ x_1=e^{-j},\ x_{10}=e^{2j}.\label{eq23}
\end{equation}
It has a closed form solution: the points are uniformly distributed on $\mathbb{S}$, with $x_n = e^{j\frac{n-4}{3}}$; that is, the angles $\arg(x_n)$ linearly interpolate between $-1$ and $2$. The proposed convex relaxation is
\begin{align}
&\minimize_{\small\substack{x_n \in \mathbb{C}\;:\; n=1,\ldots,10\\ r_{n,n+1}\in \mathbb{C}\;:\; n=1,\ldots,9}} \; \sum_{n=1}^9 \big(1- \Re(r_{n,n+1})\big)\notag\\
&\qquad \ \;\mbox{s.t.}\ \ x_1=e^{-j},\ x_{10}=e^{2j},\  \mbox{and}\label{eq24}\\ 
&\qquad \ \left[\begin{array}{ccc}
1 &x_n^*&x_{n+1}^* \\
x_n &1&r_{n,n+1}\\
x_{n+1}&r_{n,n+1}^*&1
\end{array}\right]\succcurlyeq 0,\ \forall n=1,\ldots,9.\notag
\end{align}
We solve the problem using the Chambolle--Pock algorithm \cite{cha11a,con192}, which is similar to the algorithm shown in Section~\ref{sec4}, with additional enforcement of $x_1=e^{-j}$, $x_{10}=e^{2j}$ at every iteration. Again, it turns out that the  relaxation \eqref{eq24} is tight and we obtain the exact solution to \eqref{eq23}. The interpolated signal is shown in Fig.~4, in red, and is uniform on $\mathbb{S}$, as predicted.

On the other hand, the baseline method consists in solving the convex problem
\begin{equation}
\minimize_{x_n \in \mathbb{C}\,:\, n=1,\ldots,10}\;\sum_{n=1}^9  {\textstyle\frac{1}{2}} |x_{n+1}-x_{n}|^2\ \ \mbox{s.t.}\ \ x_1=e^{-j},\ x_{10}=e^{2j},\label{eq25}
\end{equation}
and then rescaling the obtained $x_n$ as $x_n/|x_n|$ to project them on $\mathbb{S}$. We solved \eqref{eq25} using projected gradient descent, but the problem actually has a closed form solution, too: it is linear interpolation in $\mathbb{C}$, so that $x_n=\frac{n-1}{9}e^{2j}+\frac{10-n}{9}e^{-j}$. The interpolated signal is shown in Fig.~4, in blue. As we see, the angles of the $x_n$ are not uniform in $\mathbb{S}$, so that the baseline method gives a bad approximate solution to the problem \eqref{eq23}. This illustrates that the proposed convex relaxation is much tighter than the naive relaxation, which consists in reasoning in the disk $\mathbb{D}$ instead of the circle $\mathbb{S}$.

\section{Conclusion}

We proposed a new approach to smoothen or interpolate signals defined on the nonconvex complex circle, with a nonconvex formulation translating Bayesian estimation with von Mises priors, and a convex relaxation based on semidefinite programming. We showed by experiments that the proposed relaxation is tight and yields the exact solution of the nonconvex problem in several cases. This opens the door to solutions of better quality for many applications involving  circular data.

\bibliographystyle{IEEEtran}
\bibliography{IEEEabrv,../biblio.bib}

\end{document}